\documentclass[11pt]{article}
\usepackage{amssymb,latexsym}
\usepackage{epsfig}
\usepackage{eufrak}
\usepackage{amsmath}
\usepackage{mathrsfs}
\usepackage{color}

\newtheorem{theorem}{Theorem}
\newtheorem{lemma}{Lemma}

\newcommand{\be}{\begin{equation}}
\newcommand{\ee}{\end{equation}}
\newcommand{\bee}{\begin{eqnarray*}}
\newcommand{\eee}{\end{eqnarray*}}
\newcommand{\bel}{\begin{eqnarray}}
\newcommand{\eel}{\end{eqnarray}}
\newcommand{\bec}{\begin{cases}}
\newcommand{\eec}{\end{cases}}
\newcommand{\bem}{\begin{bmatrix}}
\newcommand{\eem}{\end{bmatrix}}

\newcommand{\la}{\label}
\newcommand{\li}{\left}
\newcommand{\ri}{\right}

\newcommand{\ovl}{\overline}
\newcommand{\udl}{\underline}

\newcommand{\lc}{\lceil}
\newcommand{\rc}{\rceil}
\newcommand{\lf}{\lfloor}
\newcommand{\rf}{\rfloor}

\newcommand{\vep}{\varepsilon}
\newcommand{\lm}{\lambda}

\newcommand{\de}{\delta}

\newcommand{\ga}{\gamma}
\newcommand{\Ga}{\Gamma}

\newcommand{\se}{\theta}

\newcommand{\ze}{\zeta}
\newcommand{\al}{\alpha}

\newcommand{\ro}{\rho}

\newcommand{\om}{\omega}
\newcommand{\Om}{\Omega}

\newcommand{\f}{\frac}

\newcommand{\cd}{\cdots}

\newcommand{\qqu}{\qquad}
\newcommand{\fa}{\forall}

\newcommand{\mscr}{\mathscr}

\newcommand{\mbf}{\mathbf}
\newcommand{\bb}{\mathbb}

\newcommand{\wh}{\widehat}

\newcommand{\mrm}{\mathrm}
\newcommand{\bs}{\boldsymbol}

\newcommand{\sh}{\slash}

\newcommand{\tx}{\text}

\newcommand{\iy}{\infty}

\newcommand{\bed}{\begin{description}}
\newcommand{\eed}{\end{description}}
\newcommand{\bei}{\begin{itemize}}
\newcommand{\eei}{\end{itemize}}
\newcommand{\ben}{\begin{enumerate}}
\newcommand{\een}{\end{enumerate}}
\newcommand{\bib}{\bibitem}
\newcommand{\beL}{\begin{lemma}}
\newcommand{\eeL}{\end{lemma}}
\newcommand{\beT}{\begin{theorem}}
\newcommand{\eeT}{\end{theorem}}
\newcommand{\sect}{\section}

\newcommand{\bpf}{\begin{pf}}
\newcommand{\epf}{\end{pf}}
\newcommand{\bsk}{\bigskip}
\newcommand{\bi}{\binom}

\newcommand{\pfbox}{\hfill\mbox{$\Box$}}

\newenvironment{pf}{\paragraph*{Proof{\rm.}}}{\pfbox\bigskip}

\begin{document}

\title{{\bf Estimating the Parameters of Binomial and Poisson Distributions via Multistage Sampling}
\thanks{The author had been previously working with Louisiana State University at Baton Rouge, LA 70803, USA,
and is now with Department of Electrical Engineering, Southern
University and A\&M College, Baton Rouge, LA 70813, USA; Email:
chenxinjia@gmail.com}}

\author{Xinjia Chen}

\date{September 2008}

\maketitle

\begin{abstract}

In this paper, we have developed a new class of sampling schemes for
estimating parameters of binomial and Poisson distributions. Without
any information of the unknown parameters, our sampling schemes
rigorously guarantee prescribed levels of precision and confidence.

\end{abstract}

\sect{Introduction}

The binomial and Poisson distributions are extremely useful in
numerous fields of sciences and engineering. The binomial
distribution arises in many different contexts whenever a random
variable can be hypothesized to have arisen as the number of
occurrences of a certain characteristics or property of interest in
a series of independent trials of the random phenomenon. It has been
utilized for statistical inferences about dichotomous data for more
than $250$ years.  The Poisson distribution has found an extensive
application for a wide variety of phenomena dealing with counts of
rare events (see, e.g., \cite{John, Sahai, Sahai2} and the
references therein).

The estimation of the parameters of binomial and Poisson
distributions is of practical importance and has been persistent
issues of research in statistics and other relevant fields.  Despite
the richness of literature devoted to such issues, existing
approaches suffer from the drawbacks of lacking either efficiency or
rigorousness. Such drawbacks are due to conservative bounding or
asymptotic approximation involved in the design of sampling schemes
(see, e.g., \cite{Gosh} and the references therein). To overcome the
limitations of existing methods of estimating the parameters of
binomial and Poisson distributions, we would like to propose a new
classes of multistage sampling schemes. In contrast to existing
methods, our sampling schemes require no information of the unknown
parameters and rigorously guarantee prescribed levels of precision
and confidence.

The remainder of the paper is organized as follows. In Section 2, we
present our multistage sampling schemes for estimating binomial
parameters under different precision requirements.  Section 3 is
devoted to the estimation of Poisson parameters.   Section 4 is the
conclusion.  The proofs of all theorems are given in Appendices.

Throughout this paper, we shall use the following notations. The set
of integers is denoted by $\bb{Z}$. The ceiling function and floor
function are denoted respectively by $\lc . \rc$ and $\lf . \rf$
(i.e., $\lc x \rc$ represents the smallest integer no less than $x$;
$\lf x \rf$ represents the largest integer no greater than $x$). The
gamma function is denoted by $\Ga(.)$. For any integer $m$, the
combinatoric function $\bi{m}{z}$ with respect to integer $z$ takes
value {\small $\f{ \Ga( m + 1) } { \Ga( z + 1) \Ga (m- z + 1) }$}
for $z \leq m$ and value $0$ otherwise. We use the notation $\Pr \{
. \mid \se \}$ to indicate that the associated random samples $X_1,
X_2, \cd$ are parameterized by $\se$. The parameter $\se$ in $\Pr \{
. \mid \se \}$  may be dropped whenever this can be done without
introducing confusion. The other notations will be made clear as we
proceed.

\sect{Estimation of Binomial Parameter}

Let $X$ be a Bernoulli random variable with distribution $\Pr \{ X =
1 \} = 1 - \Pr \{ X = 0 \} = p \in (0, 1)$.  It is a frequent
problem to estimate $p$ based on i.i.d. random samples $X_1, X_2,
\cd $ of $X$.  In this regard, we have developed various sampling
schemes by virtue of the following function:
\[
S_{\mrm{B}} (k, l, n, p) = \bec \sum_{i = k}^l \bi{n}{i} p^i (1 -
p)^{n - i} & \tx{for} \; p \in (0, 1),\\
0 & \tx{for} \; p \notin (0, 1).  \eec
\]

\subsection{Control of Absolute Error}

To construct an estimator satisfying an absolute error criterion
with a prescribed confidence level, we have

\beT \la{coverage_abs} Let $0 < \vep < \f{1}{2}, \; 0 < \de < 1, \;
\ze > 0$ and $\ro > 0$. Let $n_1 < n_2 < \cd < n_s$ be the ascending
arrangement of all distinct elements of {\small $ \li \{ \li \lc \li
( \f{ 2 \vep^2} { \ln \f{1}{1- \vep} } \ri )^{ 1 - \f{i}{\tau} } \f{
\ln \f{1}{\ze \de} } { 2 \vep^2 } \ri \rc : i = 0, 1, \cd, \tau \ri
\}$} with {\small  $\tau = \li \lc \f{ \ln \li ( \f{1}{ 2 \vep^2}
\ln \f{1}{1- \vep} \ri ) } { \ln (1 + \ro)} \ri \rc$. } Define
$K_\ell = \sum_{i=1}^{n_\ell} X_i$ and $\wh{\bs{p}}_\ell =
\f{K_\ell}{n_\ell}$
 for $\ell = 1, \cd, s$. Suppose the stopping
rule is that sampling is continued until $S_{\mrm{B}} (K_\ell,
n_\ell, n_\ell, \wh{\bs{p}}_\ell - \vep ) \leq \ze \de$ and
$S_{\mrm{B}} (0, K_\ell, n_\ell, \wh{\bs{p}}_\ell + \vep ) \leq \ze
\de$ for some $\ell \in \{1, \cd, s\}$. Let $\bs{\wh{p}} =
\f{\sum_{i=1}^{\mathbf{n}} X_i}{\mathbf{n}}$ where $\mathbf{n}$ is
the sample size when the sampling is terminated. Then, $ \Pr \li \{
\li | \bs{\wh{p}} - p
 \ri | < \vep \mid p \ri \} \geq 1 - \de$ for any $p \in (0, 1)$ provided that $0 < \ze \leq \f{1}{2 (\tau + 1) }$.
 \eeT

\bsk

We would like to note that if we define {\small  \[ \mscr{Q}^+ =
\bigcup_{\ell = 1}^s \li \{ \f{k}{n_\ell} + \vep \in \li (0,
\f{1}{2} \ri ) : k \in \bb{Z} \ri \} \bigcup \li \{ \f{1}{2} \ri \},
\qqu \mscr{Q}^- = \bigcup_{\ell = 1}^s \li \{ \f{k}{n_\ell} - \vep
\in \li (0, \f{1}{2} \ri ) : k \in \bb{Z} \ri \} \bigcup \li \{
\f{1}{2} \ri \}.
\]}
and decision variables $\bs{D}_\ell$ such that $\bs{D}_\ell = 1$ if
\[ S_{\mrm{B}} (K_\ell, n_\ell, n_\ell, \wh{\bs{p}}_\ell - \vep )
\leq \ze \de, \qqu S_{\mrm{B}} (0, K_\ell, n_\ell, \wh{\bs{p}}_\ell
+ \vep ) \leq \ze \de; \]
 and $\bs{D}_\ell = 0$ otherwise, then a sufficient condition to guarantee $\Pr \li \{ \li |
\bs{\wh{p}} - p \ri | < \vep \mid p \ri \}
> 1 - \de$ for any $p \in (0, 1)$ is that {\small  \bel &  &
\sum_{\ell = 1}^s \Pr \{ \wh{\bs{p}}_\ell \geq p + \vep, \;
\bs{D}_{\ell - 1} = 0, \; \bs{D}_\ell = 1 \mid p \} < \f{\de}{2}
\qqu \fa p \in
\mscr{Q}^-, \la{2D1}\\
&  & \sum_{\ell = 1}^s \Pr \{ \wh{\bs{p}}_\ell  \leq p - \vep, \;
\bs{D}_{\ell - 1} = 0, \; \bs{D}_\ell = 1 \mid p \} < \f{\de}{2}
\qqu \fa p \in \mscr{Q}^+ \la{2D2}
 \eel}
where both (\ref{2D1}) and (\ref{2D2}) are satisfied if $0 < \ze <
\f{1}{2(\tau + 1)}$.  Here we have used the double-decision-variable
method of \cite{Chen_EST}.  To determine a $\ze$ as large as
possible and thus make the sampling scheme most efficient, the
computational techniques such as bisection confidence tuning, domain
truncation, triangular partition developed in \cite{Chen_EST, ChenT}
can be applied.

\subsection{Control of Relative Error}

To construct an estimator satisfying a relative error criterion with
a prescribed confidence level, we have

\beT Let $0 < \vep < 1, \; 0 < \de < 1, \; \ze > 0$ and $\ro > 0$.
Define $\nu = \f{ \vep } { (1 + \vep) \ln (1 + \vep) - \vep}$ and
{\small $\tau = \li \lc \f{ \ln ( 1 + \nu ) }{ \ln (1 + \ro) } \ri
\rc$}. Let $\ga_1 < \ga_2 < \cd < \ga_s$ be the ascending
arrangement of all distinct elements of {\small $\li \{ \li \lc  ( 1
+ \nu )^{ \f{i}{\tau} } \f{ \ln \f{1} {\ze \de}} { \ln (1 + \vep) }
\ri \rc: i = 0, 1, \cd, \tau \ri \}$}.  For $\ell = 1, \cd, s$, let
$\wh{\bs{p}}_\ell = \f{\ga_\ell}{ \mathbf{n}_\ell }$ where
$\mathbf{n}_\ell$ is the minimum number of samples such that
$\sum_{i = 1}^{\mathbf{n}_\ell} X_i = \ga_\ell$.  Suppose the
stopping rule is that sampling is continued until {\small
$S_{\mrm{B}} \li (\ga_\ell, \mathbf{n}_\ell, \mathbf{n}_\ell,
\f{\wh{\bs{p}}_\ell}{1 + \vep} \ri ) \leq \ze \de$} and {\small
$S_{\mrm{B}} \li (0, \ga_\ell, \mathbf{n}_\ell,
\f{\wh{\bs{p}}_\ell}{1 - \vep} \ri ) \leq \ze \de$} for some $\ell
\in \{1, \cd, s\}$. Define estimator $\wh{\bs{p}} =
\f{\sum_{i=1}^{\mathbf{n}} X_i}{\mathbf{n}}$ where $\mathbf{n}$ is
the sample size when the sampling is terminated. Then, {\small $\Pr
\li \{ \li |  \f{ \wh{\bs{p}} - p } { p }  \ri | < \vep \mid p \ri
\} \geq 1 - \de$ } for any $p \in (0, 1)$  provided that $0 < \ze
\leq \f{1}{2 (\tau + 1) }$.
 \eeT

In this section, we have proposed a multistage inverse sampling plan
for estimating a binomial parameter, $p$, with relative precision.
In some situations, the cost of sampling operation may be high since
samples are obtained one by one when inverse sampling is involved.
In view of this fact, it is desirable to develop multistage
estimation methods without using inverse sampling.  For this
purpose, we have

\beT \la{noinverse} Let $0 < \vep < 1, \; 0 < \de < 1$ and $\ze
> 0 $.  Let $\tau$ be a positive integer. For $\ell = 1, 2, \cd$, define $K_\ell = \sum_{i =
1}^{n_\ell} X_i, \; \wh{\bs{p}}_\ell = \f{K_\ell} {n_\ell }$, where
$n_\ell$ is deterministic and stands for the sample size at the
$\ell$-th stage.  Suppose the stopping rule is that sampling is
continued until {\small $S_{\mrm{B}} \li (K_\ell, n_\ell, n_\ell,
\f{\wh{\bs{p}}_\ell}{1 + \vep} \ri ) \leq \ze \de_\ell$} and {\small
$S_{\mrm{B}} \li (0, K_\ell, n_\ell, \f{\wh{\bs{p}}_\ell}{1 - \vep}
\ri ) \leq \ze \de_\ell$} for some $\ell$, where $\de_\ell = \de$
for $1 \leq \ell \leq \tau$ and $\de_\ell = \de 2^{\tau - \ell}$ for
$\ell > \tau$. Define estimator $\wh{\bs{p}} =
\wh{\bs{p}}_{\bs{l}}$, where $\bs{l}$ is the index of stage at which
the sampling is terminated. Then, $\Pr \{ \bs{l} < \iy \} = 1$ and
{\small $\Pr \li \{ \li | \f{ \wh{\bs{p}} - p } { p } \ri | \leq
\vep \mid p \ri \} \geq 1 - \de$} for any $p \in (0, 1)$ provided
that $2 (\tau + 1 ) \ze \leq 1$ and $\inf_{\ell > 0} \f{n_{\ell +
1}}{n_\ell} > 0$.

\eeT

\subsection{Control of Absolute and Relative Errors}

To construct an estimator satisfying a mixed criterion in terms of
absolute and relative errors with a prescribed confidence level, we
have

\beT \la{coverage_mixed}  Let $0 < \de < 1, \; \ze > 0$ and $\ro >
0$. Let $\vep_a$ and $\vep_r$ be positive numbers such that $0 <
\vep_a < \f{35}{94}$ and {\small $\f{70 \vep_a}{35 - 24 \vep_a } <
\vep_r < 1$}. Define {\small $\nu = \f{ \vep_a + \vep_r \vep_a -
\vep_r} { \vep_r \ln (1 + \vep_r) } \ln \li ( 1 + \f{\vep_r^2} {
\vep_r - \vep_a - \vep_r \vep_a } \ri )$} and {\small $\tau = \li
\lf \f{\ln (1 + \nu)}{ \ln (1 + \ro) } \ri \rf$}. Let $n_1 < n_2 <
\cd < n_s$ be the ascending arrangement of all distinct elements of
{\small $\li \{ \li \lc ( 1 + \nu )^{ \f{i} { \tau } } \f{ \ln \f{1}
{ \ze \de } }{\ln (1 + \vep_r)} \ri \rc: \tau \leq i \leq 0 \ri
\}$}. Define {\small $K_\ell = \sum_{i=1}^{n_\ell} X_i, \;
\wh{\bs{p}}_\ell = \f{K_\ell}{n_\ell}, \; \udl{\bs{p}}_\ell = \min
\{ \wh{\bs{p}}_\ell - \vep_a, \; \f{\wh{\bs{p}}_\ell}{1 + \vep_r}
\}$} and {\small $\ovl{\bs{p}}_\ell = \max \{ \wh{\bs{p}}_\ell +
\vep_a, \; \f{\wh{\bs{p}}_\ell}{1 - \vep_r} \}$}
 for $\ell = 1, \cd, s$.  Suppose the stopping rule is that sampling is continued until
 {\small $S_{\mrm{B}} (K_\ell, n_\ell, n_\ell,  \udl{\bs{p}}_\ell )  \leq \ze \de$} and {\small
$S_{\mrm{B}} (0, K_\ell, n_\ell,  \ovl{\bs{p}}_\ell )  \leq \ze
\de$} for some $\ell \in \{1, \cd, s\}$.   Let {\small $\wh{\bs{p}}
= \f{\sum_{i=1}^{\mathbf{n}} X_i} {\mathbf{n}}$} where $\mathbf{n}$
is the sample size when the sampling is terminated. Then, {\small
$\Pr \li \{ \li | \wh{\bs{p}} - p
 \ri | < \vep_a  \; \mrm{or} \; \li |
\f{\wh{\bs{p}} - p } {p } \ri | < \vep_r  \mid p \ri \} \geq 1 -
\de$
 }
for any $p \in (0, 1)$ provided that $0 < \ze \leq \f{1}{2 (1 -
\tau)}$. \eeT

\sect{Estimation of Poisson Parameter}

Let $X$ be a Poisson random variable with mean value $\lm
> 0$.  It is an important problem to estimate $\lm$ based on i.i.d.
random samples $X_1, X_2, \cd $ of $X$. In this regard, we have
developed a sampling scheme by virtue of the following function:
\[
S_{\mrm{P}} (k, l, n, \lm) = \bec \sum_{i = k}^l \f{ (n \lm)^i e^{- n \lm}} { i! }  & \tx{for} \; \lm > 0,\\
0 & \tx{for} \; \lm \leq 0.  \eec
\]
As can be seen at below, our sampling scheme produces an estimator
satisfying a mixed criterion in terms of absolute and relative
errors with a prescribed confidence level.

\beT Let $0 < \vep_a < 1, \; 0 < \vep_r < 1, \; 0 < \de < 1, \; \ze
> 0$ and $\ro > 0$.  Let $n_1 < n_2 < \cd < n_s$ be the ascending arrangement of
all distinct elements of {\small $\li \{ \li \lc \nu^{ \f{i} { \tau
} }  \ln \f{1} { \ze \de } \ri \rc: i = 0, 1, \cd, \tau \ri \}$}
with $\nu =  \f{\vep_r}{ \vep_a [( 1 + \vep_r ) \ln \li ( 1 + \vep_r
\ri ) - \vep_r ] } $ and {\small $\tau  = \li \lc \f{ \ln \nu } {
\ln ( 1 + \ro) } \ri \rc$}.   Define {\small $K_\ell =
\sum_{i=1}^{n_\ell} X_i, \; \wh{\bs{p}}_\ell = \f{K_\ell}{n_\ell},
\; \ovl{\bs{\lm}}_\ell = \max \li \{ \wh{\bs{\lm}}_\ell + \vep_a, \;
\f{\wh{\bs{\lm}}_\ell}{1 - \vep_r} \ri \}$} and {\small
$\udl{\bs{\lm}}_\ell = \min \{ \wh{\bs{\lm}}_\ell - \vep_a, \;
\f{\wh{\bs{\lm}}_\ell}{1 + \vep_r} \}$} for $\ell \in \{1, \cd,
s\}$.  Suppose the stopping rule is that sampling is continued until
{\small $S_{\mrm{P}} (0, K_\ell - 1, n_\ell, \udl{\bs{\lm}}_\ell)
\geq 1 - \ze \de$} and {\small $S_{\mrm{P}} (0, K_\ell, n_\ell,
\ovl{\bs{\lm}}_\ell) \leq \ze \de$} for some $\ell \in \{1, \cd,
s\}$.  Let {\small $\wh{\bs{\lm}} = \f{ \sum_{i=1}^{\mathbf{n}} X_i
} {\mathbf{n}}$} where $\mathbf{n}$ is the sample size when the
sampling is terminated. Then, {\small $\Pr \{  | \wh{\bs{\lm}} - \lm
 | < \vep_a  \; \mrm{or} \;  | \wh{\bs{\lm}} - \lm
 | < \vep_r \lm \mid \lm \} \geq 1 - \de$}
 for any $\lm \in (0, \iy)$ provided that $0 < \ze \leq \f{1}{2 (\tau + 1) }$.  \eeT

For the purpose of estimating Poisson parameter, $\lm$, with an
absolute precision, we have

 \beT \la{noinversePabs} Let $\vep > 0, \; 0 < \de < 1$ and $\ze
> 0 $.  Let $\tau$ be a positive integer. For $\ell = 1, 2, \cd$, define $K_\ell = \sum_{i =
1}^{n_\ell} X_i, \; \wh{\bs{\lm}}_\ell = \f{K_\ell} {n_\ell }$,
where $n_\ell$ is deterministic and stands for the sample size at
the $\ell$-th stage.  Suppose the stopping rule is that sampling is
continued until {\small $S_{\mrm{P}} \li (0, K_\ell - 1, n_\ell,
\wh{\bs{\lm}}_\ell -  \vep \ri ) \geq 1 - \ze \de_\ell$} and {\small
$S_{\mrm{P}} \li (0, K_\ell, n_\ell, \wh{\bs{\lm}}_\ell + \vep \ri )
\leq \ze \de_\ell$} for some $\ell$, where $\de_\ell = \de$ for $1
\leq \ell \leq \tau$ and $\de_\ell = \de 2^{\tau - \ell}$ for $\ell
> \tau$. Define estimator $\wh{\bs{\lm}} = \wh{\bs{\lm}}_{\bs{l}}$,
where $\bs{l}$ is the index of stage at which the sampling is
terminated. Then, $\Pr \{ \bs{l} < \iy \} = 1$ and {\small $\Pr \li
\{ \li | \wh{\bs{\lm}} - \lm \ri | \leq \vep \mid \lm \ri \} \geq 1
- \de$} for any $\lm \in (0, \iy)$ provided that $2 (\tau + 1 ) \ze
\leq 1$ and $\inf_{\ell > 0} \f{n_{\ell + 1}}{n_\ell} > 0$.

\eeT

For the purpose of estimating Poisson parameter, $\lm$, with a
relative precision, we have

 \beT \la{noinversePrev} Let $0 < \vep < 1, \; 0 < \de < 1$ and $\ze
> 0 $.  Let $\tau$ be a positive integer. For $\ell = 1, 2, \cd$, define $K_\ell = \sum_{i =
1}^{n_\ell} X_i, \; \wh{\bs{\lm}}_\ell = \f{K_\ell} {n_\ell }$,
where $n_\ell$ is deterministic and stands for the sample size at
the $\ell$-th stage.  Suppose the stopping rule is that sampling is
continued until {\small $S_{\mrm{P}} \li (0, K_\ell - 1, n_\ell,
\f{\wh{\bs{\lm}}_\ell}{1 + \vep} \ri ) \geq 1 - \ze \de_\ell$} and
{\small $S_{\mrm{P}} \li (0, K_\ell, n_\ell,
\f{\wh{\bs{\lm}}_\ell}{1 - \vep} \ri ) \leq \ze \de_\ell$} for some
$\ell$, where $\de_\ell = \de$ for $1 \leq \ell \leq \tau$ and
$\de_\ell = \de 2^{\tau - \ell}$ for $\ell > \tau$. Define estimator
$\wh{\bs{\lm}} = \wh{\bs{\lm}}_{\bs{l}}$, where $\bs{l}$ is the
index of stage at which the sampling is terminated. Then, $\Pr \{
\bs{l} < \iy \} = 1$ and {\small $\Pr \li \{ \li | \f{ \wh{\bs{\lm}}
- \lm } { \lm } \ri | \leq \vep \mid \lm \ri \} \geq 1 - \de$} for
any $\lm \in (0, \iy)$ provided that $2 (\tau + 1 ) \ze \leq 1$ and
$\inf_{\ell > 0} \f{n_{\ell + 1}}{n_\ell} > 0$.

\eeT

 \bsk

 Again, as we mentioned after the presentation of Theorem 1,
 we would like to note that the computational techniques such as the double-decision-variable method, bisection confidence tuning,
domain truncation, triangular partition developed in
 \cite{Chen_EST, ChenT} can be applied to reduce the conservatism of the sampling
 schemes described by Theorems 2 to 7.

 With regard to the tightness of the double-decision-variable method,
we can develop results similar to Theorems 8, 13, 18 and 22 of
\cite{Chen_EST}.

With regard to the asymptotic performance of our sampling schemes,
we can develop results similar to Theorems 9, 14, 19 and 23 of
\cite{Chen_EST}.

 \sect{Conclusion}

In this paper, we have developed new multistage sampling schemes for
estimating the parameters of binomial and Poisson distributions. Our
new methods rigorously guarantee prescribed levels of precision and
confidence.

\appendix

\sect{Proof of Theorem 1}

In the course of proof, we need to use function
\[
\mscr{M}_{\mrm{B}} (z,\mu) = \bec z \ln \f{\mu}{z} + (1 - z) \ln
\f{1 - \mu}{1 - z} &
\tx{for} \; z \in (0,1) \; \tx{and} \; \mu \in (0, 1),\\
\ln(1-\mu) & \tx{for} \; z = 0 \; \tx{and} \; \mu \in (0, 1),\\
\ln \mu &  \tx{for} \; z = 1 \; \tx{and} \; \mu \in (0, 1),\\
- \iy &  \tx{for} \; z \in [0, 1] \; \tx{and} \; \mu \notin (0, 1).
\eec
\]

We need some preliminary results. The following classical result is
due to Hoeffding \cite{Hoeffding}.

\beL  \la{Hoe} Let $\ovl{X}_n = \f{\sum_{i=1}^n X_i}{n}$ where $X_1,
\; \cd,\; X_n$ are i.i.d.  random variables such that $0 \leq X_i
\leq 1$ and $\bb{E}[X_i] = \mu \in (0, 1)$ for $i = 1, \; \cd, n$.
Then, $\Pr \li \{ \ovl{X}_n \geq z \ri \} \leq \exp \li (n
\mscr{M}_{\mrm{B}} \li (z, \mu \ri ) \ri )$ for any $z \in (\mu,
1)$. Similarly, $\Pr \li \{ \ovl{X}_n \leq z \ri \} \leq \exp \li (n
\mscr{M}_{\mrm{B}} \li (z, \mu \ri ) \ri )$ for any $z \in (0,
\mu)$. \eeL

The following lemma can be readily derived from Lemma \ref{Hoe}.

\beL \la{decb}

{\small $S_{\mrm{B}} (0, k, n, p) \leq \exp ( n \mscr{M}_{\mrm{B}} (
\f{k}{n}, p  )  )$} for $0 \leq k \leq n p$. Similarly, {\small
$S_{\mrm{B}} (k, n, n, p) \leq \exp  ( n \mscr{M}_{\mrm{B}} (
\f{k}{n}, p  )  )$} for $n p \leq k \leq n$. \eeL

\beL

\la{absineq1}

Let $K = \sum_{i=1}^n X_i$ where $X_1, \cd, X_n$ are i.i.d.
Bernoulli random variables such that $\Pr \{ X_i = 1 \} = 1 - \Pr \{
X_i = 0 \} = p \in (0, 1)$ for $i = 1, \cd, n$. Then, $\Pr \li \{
S_{\mrm{B}} \li ( 0, K, n, p \ri ) \leq \al \ri \} \leq \al$ for any
$\al > 0$. \eeL

\bpf If $\{S_{\mrm{B}}  (0, K, n, p ) \leq \al \}$ is an impossible
event, then $\Pr  \{ S_{\mrm{B}} (0, K, n, p ) \leq \al \} = 0 <
\al$. Otherwise, if $ \{S_{\mrm{B}} \li (0, K, n, p \ri ) \leq \al
 \}$ is a possible event, then there exists an integer $k^* = \max \{k : 0 \leq k \leq n, \;
S_{\mrm{B}}  ( 0, k, n, p ) \leq \al \}$
 and it follows that $ \Pr \{ S_{\mrm{B}}  (0, K, n, p ) \leq \al \} = S_{\mrm{B}} (0,
k^*, n, p  ) \leq \al$. The proof is thus completed. \epf

\beL \la{absineq2}
 Let $K = \sum_{i=1}^n X_i$ where $X_1, \cd, X_n$
are i.i.d. Bernoulli random variables such that $\Pr \{ X_i = 1 \} =
1 - \Pr \{ X_i = 0 \} = p \in (0, 1)$ for $i = 1, \cd, n$.   Then,
$\Pr \li \{ S_{\mrm{B}} \li ( K, n, n, p \ri ) \leq \al \ri \} \leq
\al$ for any $\al > 0$.

\eeL

\bpf If $\{S_{\mrm{B}} (K, n, n, p ) \leq \al \}$ is an impossible
event, then $\Pr  \{ S_{\mrm{B}}  (  K, n, n, p ) \leq \al \} = 0 <
\al$. Otherwise, if $\{S_{\mrm{B}} (K, n, n, p ) \leq \al \}$ is a
possible event, then there exists an integer {\small $k_\star = \min
\{k : 0 \leq k \leq n, \; S_{\mrm{B}} ( k, n, n, p ) \leq \al \}$}
and it follows that $\Pr \{ S_{\mrm{B}} ( K, n, n, p ) \leq \al \} =
S_{\mrm{B}} (k_\star, n, n, p ) \leq \al$. The proof is thus
completed. \epf

\beL

\la{lemax}

Both $\mscr{M}_{\mrm{B}} (z, z - \vep)$ and $\mscr{M}_{\mrm{B}} (z,
z + \vep)$  are no greater than $- 2 \vep^2$ for $0 \leq z \leq 1$.

\eeL

\beL \la{absDS1}

$\Pr \li \{ S_{\mrm{B}} \li (K_s, n_s, n_s, \wh{\bs{p}}_s - \vep \ri
) \leq \ze \de \ri \} = \Pr \li \{ S_{\mrm{B}} \li (0, K_s, n_s,
\wh{\bs{p}}_s + \vep \ri ) \leq \ze \de \ri \} = 1$.

 \eeL

 \bpf
By the definition of sample sizes, we have $n_s = \li \lc \f{\ln
(\ze \de) }{ - 2 \vep^2}  \ri \rc \geq \f{\ln (\ze \de) }{ - 2
\vep^2}$ and consequently  $\f{\ln (\ze \de) }{n_s} \geq - 2 \vep^2$
By Lemmas \ref{decb} and \ref{lemax}, we have {\small \[ \Pr \li \{
S_{\mrm{B}} \li (K_s, n_s, n_s, \wh{\bs{p}}_s - \vep \ri ) \leq \ze
\de \ri \} \geq \Pr \li \{ \mscr{M}_{\mrm{B}} \li (\wh{\bs{p}}_s,
\wh{\bs{p}}_s - \vep \ri ) \leq \f{\ln (\ze \de)}{n_s} \ri \} \geq
\Pr \li \{ \mscr{M}_{\mrm{B}} \li (\wh{\bs{p}}_s, \wh{\bs{p}}_s -
\vep \ri ) \leq - 2 \vep^2 \ri \} = 1, \]
\[ \Pr \li \{
S_{\mrm{B}} \li (0, K_s, n_s, \wh{\bs{p}}_s + \vep \ri ) \leq \ze
\de \ri \} \geq \Pr \li \{ \mscr{M}_{\mrm{B}} \li (\wh{\bs{p}}_s,
\wh{\bs{p}}_s + \vep \ri ) \leq \f{\ln (\ze \de)}{n_s} \ri \} \geq
\Pr \li \{ \mscr{M}_{\mrm{B}} \li (\wh{\bs{p}}_s, \wh{\bs{p}}_s +
\vep \ri ) \leq - 2 \vep^2 \ri \} = 1
\]}
which immediately implies the lemma.

 \epf

 \beL
 \la{good1}
$\Pr \{ p \leq \wh{\bs{p}}_\ell - \vep, \;  S_{\mrm{B}} \li (K_\ell,
n_\ell, n_\ell, \wh{\bs{p}}_\ell - \vep \ri ) \leq \ze \de \} \leq
\ze \de$ for $\ell = 1, \cd, s$.  \eeL

\bpf

Since $S_{\mrm{B}} (k, n, n, p)$ is monotonically increasing with
respect to $p \in (0, 1)$, we have $\{ p \leq \wh{\bs{p}}_\ell -
\vep, \;  S_{\mrm{B}} \li (K_\ell, n_\ell, n_\ell, \wh{\bs{p}}_\ell
- \vep \ri ) \leq \ze \de \} \subseteq \{  S_{\mrm{B}} \li (K_\ell,
n_\ell, n_\ell, p \ri ) \leq \ze \de \}$.  Hence, by Lemma
\ref{absineq2}, we have
\[
\Pr \li \{ p \leq \wh{\bs{p}}_\ell - \vep, \;  S_{\mrm{B}} \li
(K_\ell, n_\ell, n_\ell, \wh{\bs{p}}_\ell - \vep \ri ) \leq \ze \de
\ri \} \leq \Pr \li \{  S_{\mrm{B}} \li (K_\ell, n_\ell, n_\ell, p
\ri ) \leq \ze \de \ri \} \leq \ze \de
\]
for $\ell = 1, \cd, s$. \epf

 \beL
 \la{good2}
$\Pr \{ p \geq \wh{\bs{p}}_\ell + \vep, \;  S_{\mrm{B}} \li (0,
K_\ell, n_\ell, \wh{\bs{p}}_\ell + \vep \ri ) \leq \ze \de \} \leq
\ze \de$ for $\ell = 1, \cd, s$.  \eeL

\bpf

Since $S_{\mrm{B}} (0, k, n, p)$ is monotonically decreasing with
respect to $p \in (0, 1)$, we have $\{ p \geq \wh{\bs{p}}_\ell +
\vep, \;  S_{\mrm{B}} \li (0, K_\ell, n_\ell, \wh{\bs{p}}_\ell +
\vep \ri ) \leq \ze \de \} \subseteq \{ S_{\mrm{B}} \li (0, K_\ell,
n_\ell, p \ri ) \leq \ze \de \}$.  Hence, by Lemma \ref{absineq1},
we have
\[
\Pr \li \{ p \geq \wh{\bs{p}}_\ell + \vep, \;  S_{\mrm{B}} \li (0,
K_\ell, n_\ell, \wh{\bs{p}}_\ell + \vep \ri ) \leq \ze \de \ri \}
\leq \Pr \li \{  S_{\mrm{B}} \li (0, K_\ell, n_\ell, p \ri ) \leq
\ze \de \ri \} \leq \ze \de
\]
for $\ell = 1, \cd, s$. \epf

\bsk

Now we are in a position to prove Theorem 1.  As a direct
consequence of $\vep \in \li (0, \f{1}{2} \ri )$, we have  $\ln
\f{1}{1- \vep}
> 2 \vep^2$ and thus $\tau  \geq 1$.  This shows that the sample sizes $n_1,
\cd, n_s$ are well-defined.  By Lemma \ref{absDS1}, the sampling
must stop at some stage with index $\ell \in \{1, \cd, s\}$.
Therefore, the sampling scheme is well-defined.  By Lemmas
\ref{good1}, \ref{good2} and the definition of the stopping rule, we
have \bee \Pr \{  | \wh{\bs{p}} - p | \geq \vep \} & = & \Pr \{ p
\leq \wh{\bs{p}} - \vep \} +  \Pr \{ p \geq
\wh{\bs{p}} + \vep \}\\
& \leq & \sum_{\ell = 1}^s \Pr \li \{ p \leq \wh{\bs{p}}_\ell -
\vep, \; S_{\mrm{B}} \li (K_\ell, n_\ell, n_\ell, \wh{\bs{p}}_\ell -
\vep \ri ) \leq \ze \de \ri \} \\
&  & +  \sum_{\ell = 1}^s \Pr \li \{ p \geq \wh{\bs{p}}_\ell + \vep,
\; S_{\mrm{B}} \li (0, K_\ell, n_\ell, \wh{\bs{p}}_\ell + \vep \ri )
\leq \ze \de \ri \}\\
& \leq &  s \ze \de + s \ze \de = 2 s \ze \de \leq 2 (\tau + 1) \de,
\eee from which it can be seen that $\Pr \{ | \wh{\bs{p}} - p | <
\vep \}
> 1 - \de$ if $0 < \ze < \f{1}{2 (\tau + 1) }$.  This concludes the proof
of Theorem 1.

\sect{Proof of Theorem 2}

\beL

Let $\ga$ be a positive integer.  Let $\bs{n}$ be the minimum
integer such that $\sum_{i = 1}^{\bs{n}} X_i = \ga$ where $X_1, X_2,
\cd$ is a sequence of i.i.d. Bernoulli random variables such that
$\Pr \{ X_i = 1 \} = 1 - \Pr \{ X_i = 0 \} = p \in (0, 1)$ for any
positive integer $i$. Then, $\Pr \{ S_{\mrm{B}} (0, \ga, \bs{n}, p )
\leq \al \} \leq \al$ for any $\al > 0$. \eeL

\bpf  Since $\Pr \{ \bs{n} \geq m \} = S_{\mrm{B}} (0, \ga, m, p )$
and $\lim_{m \to \iy} \Pr \{ \bs{n} \geq m \} = 0$, there exists an
integer $m^* \geq r$ such that $S_{\mrm{B}} (0, \ga, m, p ) \leq
\al$ for any integer $m \geq m^*$ and that $S_{\mrm{B}} (0, \ga, m,
p ) > \al$ for $r \leq m < m^*$. Hence, $\Pr \{ S_{\mrm{B}} (0, \ga,
\bs{n}, p ) \leq \al \} = \Pr \{ \bs{n} \geq m^* \} = S_{\mrm{B}}
(0, \ga, m^*, p ) \leq \al$.

\epf

\beL

Let $\ga$ be a positive integer.  Let $\bs{n}$ be the minimum
integer such that $\sum_{i = 1}^{\bs{n}} X_i = \ga$ where $X_1, X_2,
\cd$ is a sequence of i.i.d. Bernoulli random variables  such that
$\Pr \{ X_i = 1 \} = 1 - \Pr \{ X_i = 0 \} = p \in (0, 1)$ for any
positive integer $i$. Then, $\Pr \{ S_{\mrm{B}} (\ga, \bs{n},
\bs{n}, p ) \leq \al \} \leq \al$ for any $\al > 0$. \eeL

\bpf  Note that $\Pr \{ \bs{n} \leq m \} = S_{\mrm{B}} (\ga, m,  m,
p )$.  In the case that $S_{\mrm{B}} (\ga, \ga,  \ga, p ) > \al$, we
have $S_{\mrm{B}} (\ga, m,  m, p ) \geq S_{\mrm{B}} (\ga, \ga,  \ga,
p )
> \al$ for any integer  $m \geq \ga$. Thus, $\Pr \{
S_{\mrm{B}} (\ga, \bs{n}, \bs{n}, p ) \leq \al \}  = 0 < \al$.   In
the case that $S_{\mrm{B}} (\ga, \ga, \ga, p ) \leq \al$, there
exists an integer $m^*$ such that $S_{\mrm{B}} ( \ga, m, m, p ) >
\al$ for any integer $m
> m^*$ and that $S_{\mrm{B}} (\ga, m, m, p ) \leq \al$ for
$\ga \leq m \leq m^*$.
 Hence, $\Pr \{ S_{\mrm{B}} (\ga, \bs{n}, \bs{n},  p ) \leq \al \} = \Pr \{
\bs{n} \leq m^* \} = S_{\mrm{B}} (\ga, m^*, m^*, p ) \leq \al$.

\epf

Now we need to introduce function
\[
\mscr{M}_{\mrm{I}} (z,\mu) = \bec  \ln \f{\mu}{z} + \li ( \f{1}{z} -
1 \ri ) \ln \f{1 - \mu}{1 - z} &
\tx{for} \; z \in (0,1) \; \tx{and} \; \mu \in (0, 1),\\
\ln \mu &  \tx{for} \; z = 1 \; \tx{and} \; \mu \in (0, 1),\\
- \iy & \tx{for} \; z = 0 \; \tx{and} \; \mu \in (0, 1),\\
- \iy &  \tx{for} \; z \in [0, 1] \; \tx{and} \; \mu \notin (0, 1).
\eec
\]
The following results, stated as Lemmas \ref{decr} and
\ref{compare}, have been established by Chen in \cite{Chen_EST}.

\beL \la{decr}
 Let $0 < \vep < 1$.  Then,  {\small
$\mscr{M}_{\mrm{I}} \li ( z, \f{z}{1 + \vep} \ri )$} is
monotonically decreasing with respect to $z \in (0, 1)$. \eeL

\beL \la{compare}
$\mscr{M}_{\mrm{I}} \li ( z, \f{z}{1 + \vep} \ri )
> \mscr{M}_{\mrm{I}} \li ( z, \f{z}{1 - \vep} \ri )$ for $0 < z < 1
- \vep < 1$. \eeL

\beL \la{revDS1} $\Pr \li \{ S_{\mrm{B}} \li (\ga_s, \mathbf{n}_s,
\mathbf{n}_s, \f{\wh{\bs{p}}_s}{1 + \vep} \ri ) \leq \ze \de \ri \}
= \Pr \li \{ S_{\mrm{B}} \li (0, \ga_s, \mathbf{n}_s,
\f{\wh{\bs{p}}_s}{1 - \vep} \ri ) \leq \ze \de \ri \} = 1$.

 \eeL

 \bpf

 By Lemma \ref{decb},
{\small  \bel
 \Pr \li \{ S_{\mrm{B}} \li (\ga_s, \mathbf{n}_s, \mathbf{n}_s,
\f{\wh{\bs{p}}_s}{1 + \vep} \ri ) \leq \ze \de \ri \} & \geq & \Pr
\li \{  \mathbf{n}_s \mscr{M}_{\mrm{B}} \li ( \f{ \ga_s } {
\mathbf{n}_s }, \f{\wh{\bs{p}}_s}{1 + \vep} \ri )  \leq \ln (\ze \de) \ri \} \nonumber\\
& = & \Pr \li \{  \f{\ga_s}{\wh{\bs{p}}_s} \mscr{M}_{\mrm{B}} \li (
\wh{\bs{p}}_s, \f{\wh{\bs{p}}_s}{1 + \vep} \ri ) \leq \ln (\ze \de)
\ri \} \nonumber\\
& = & \Pr \li \{  \mscr{M}_{\mrm{I}} \li ( \wh{\bs{p}}_s,
\f{\wh{\bs{p}}_s}{1 + \vep} \ri ) \leq \f{ \ln (\ze \de) }{\ga_s}
\ri \}. \la{eqrev1} \eel} Making use of Lemma \ref{decr} and the
fact {\small $\lim_{z \to 0} \mscr{M}_{\mrm{I}} \li ( z, \f{z}{1 +
\vep} \ri ) = \f{\vep}{1 + \vep} - \ln(1 + \vep)$},  we have {\small
$\mscr{M}_{\mrm{I}} \li ( z, \f{z}{1 + \vep} \ri ) < \f{\vep}{1 +
\vep} - \ln(1 + \vep)$} for any $z \in (0, 1]$. Consequently,
{\small $\li \{ \mscr{M}_{\mrm{I}} \li ( \wh{\bs{p}}_s,
\f{\wh{\bs{p}}_s}{1 + \vep} \ri ) \leq \f{\vep}{1 + \vep} - \ln(1 +
\vep) \ri \}$} is a sure event because $0 < \wh{\bs{p}}_s (\om) \leq
1$ for any $\om \in \Om$.  By the definition of $\ga_s$, we have
{\small \[ \ga_s = \li \lc \f{ \ln (\ze \de) } { \f{\vep}{1 + \vep}
- \ln(1 + \vep) } \ri \rc \geq \f{ \ln (\ze \de)  } { \f{\vep}{1 +
\vep} - \ln(1 + \vep) }.
\]}
Since $\f{\vep}{1 + \vep} - \ln(1 + \vep) < 0$ for any $\vep \in (0,
1)$, we have $\f{ \ln (\ze \de)  } { \ga_s } \geq \f{\vep}{1 + \vep}
- \ln(1 + \vep)$.  Hence, {\small \be \la{eqrev2} \Pr \li \{
 \mscr{M}_{\mrm{I}} \li ( \wh{\bs{p}}_s, \f{\wh{\bs{p}}_s}{1 + \vep}
\ri ) \leq \f{ \ln (\ze \de) }{\ga_s}  \ri \} \geq \Pr \li \{
\mscr{M}_{\mrm{I}} \li ( \wh{\bs{p}}_s, \f{\wh{\bs{p}}_s}{1 + \vep}
\ri ) \leq \f{\vep}{1 + \vep} - \ln(1 + \vep) \ri \} = 1. \ee}
Combining (\ref{eqrev1}) and (\ref{eqrev2}) yields {\small $\Pr \li
\{ S_{\mrm{B}} \li (\ga_s, \mathbf{n}_s, \mathbf{n}_s,
\f{\wh{\bs{p}}_s}{1 + \vep} \ri ) \leq \ze \de \ri \} = 1$}.

Similarly, by Lemmas \ref{decb} and \ref{compare}, {\small  \bel
 1 \geq \Pr \li \{ S_{\mrm{B}} \li (0, \ga_s, \mathbf{n}_s,
\f{\wh{\bs{p}}_s}{1 - \vep} \ri ) \leq \ze \de \ri \} & \geq & \Pr
\li \{  \mathbf{n}_s \mscr{M}_{\mrm{B}} \li ( \f{ \ga_s } {
\mathbf{n}_s }, \f{\wh{\bs{p}}_s}{1 - \vep} \ri )  \leq \ln (\ze \de) \ri \} \nonumber\\
& = & \Pr \li \{  \f{\ga_s}{\wh{\bs{p}}_s} \mscr{M}_{\mrm{B}} \li (
\wh{\bs{p}}_s, \f{\wh{\bs{p}}_s}{1 - \vep} \ri ) \leq \ln (\ze \de)
\ri \} \nonumber\\
& = & \Pr \li \{  \mscr{M}_{\mrm{I}} \li ( \wh{\bs{p}}_s,
\f{\wh{\bs{p}}_s}{1 - \vep} \ri ) \leq \f{ \ln (\ze \de)
}{\ga_s} \ri \} \nonumber\\
& \geq & \Pr \li \{  \mscr{M}_{\mrm{I}} \li ( \wh{\bs{p}}_s,
\f{\wh{\bs{p}}_s}{1 + \vep} \ri ) \leq \f{ \ln (\ze \de) }{\ga_s}
\ri \} = 1. \la{eqrev3} \eel} This completes the proof of the lemma.

 \epf

By a similar argument as that of Lemma \ref{good1}, we have Lemma
\ref{goodrev1} as follows.
  \beL
  \la{goodrev1}
{\small $\Pr \li \{ p \leq \f{\wh{\bs{p}}_\ell }{1 + \vep}, \;
S_{\mrm{B}} \li (\ga_\ell, \mbf{n}_\ell, \mbf{n}_\ell,
\f{\wh{\bs{p}}_\ell}{1 + \vep} \ri ) \leq \ze \de \ri \} \leq \ze
\de$} for $\ell = 1, \cd, s$. \eeL

By a similar argument as that of Lemma \ref{good2}, we have Lemma
\ref{goodrev2} as follows.

 \beL
  \la{goodrev2}
{\small $\Pr \li \{ p \geq \f{\wh{\bs{p}}_\ell}{1 - \vep}, \;
S_{\mrm{B}} \li (0, \ga_\ell, \mbf{n}_\ell, \f{\wh{\bs{p}}_\ell}{1 -
\vep} \ri ) \leq \ze \de \ri \} \leq \ze \de$} for $\ell = 1, \cd,
s$. \eeL

\bsk

Now we are in a position to prove Theorem 2.  Since $\ln (1 + \vep)
> \f{\vep}{1 + \vep}$ for any $\vep \in \li (0, 1 \ri )$, we have  $\nu > 0$
 and thus $\tau  \geq 1$.  This shows that the sample sizes $n_1,
\cd, n_s$ are well-defined.  By Lemma \ref{revDS1}, the sampling
must stop at some stage with index $\ell \in \{1, \cd, s\}$.
Therefore, the sampling scheme is well-defined.  By Lemmas
\ref{goodrev1}, \ref{goodrev2} and the definition of the stopping
rule, we have \bee \Pr \{  | \wh{\bs{p}} - p | \geq \vep \} & = &
\Pr \{ p \leq \wh{\bs{p}} \sh (1 + \vep) \} +  \Pr \{ p \geq
\wh{\bs{p}} \sh (1 - \vep) \}\\
& \leq & \sum_{\ell = 1}^s \Pr \li \{ p \leq \wh{\bs{p}}_\ell \sh (
1 + \vep), \; S_{\mrm{B}} \li (K_\ell, n_\ell, n_\ell,
\wh{\bs{p}}_\ell  \sh (1 + \vep) \ri ) \leq \ze \de \ri \} \\
&  & +  \sum_{\ell = 1}^s \Pr \li \{ p \geq \wh{\bs{p}}_\ell \sh (1
- \vep), \; S_{\mrm{B}} \li (0, K_\ell, n_\ell, \wh{\bs{p}}_\ell \sh
(1 - \vep ) \ri )
\leq \ze \de \ri \}\\
& \leq &  s \ze \de + s \ze \de = 2 s \ze \de \leq 2 (\tau + 1) \de,
\eee from which it can be seen that $\Pr \{ | \wh{\bs{p}} - p | <
\vep p \} > 1 - \de$ if $0 < \ze < \f{1}{2 (\tau + 1) }$.  This
concludes the proof of Theorem 2.

\sect{Proof of Theorem 4}

The following result, stated as Lemma \ref{lem88}, has been
established by Chen in \cite{Chen_EST}.

\beL \la{lem88} $\Pr \li \{ \mscr{M}_{\mrm{B}} \li ( \wh{\bs{p}}_s,
\udl{\bs{p}}_s \ri ) \leq \f{\ln (\ze \de) } { n_s }, \;
\mscr{M}_{\mrm{B}} \li ( \wh{\bs{p}}_s, \ovl{\bs{p}}_s \ri ) \leq
\f{\ln (\ze \de) } { n_s } \ri \} = 1$.  \eeL

\beL \la{mixDS1} $\Pr \{  S_{\mrm{B}} (K_s, n_s, n_s, \udl{\bs{p}}_s
) \leq \ze \de \} = \Pr \{ S_{\mrm{B}} (0, K_s, n_s, \ovl{\bs{p}}_s
) \leq \ze \de \} = 1$. \eeL

\bpf

By Lemmas \ref{decb} and \ref{lem88},
\[
1 \geq \Pr \li \{  S_{\mrm{B}} (K_s, n_s, n_s, \udl{\bs{p}}_s ) \leq
\ze \de \ri \} \geq \Pr \li \{  n_s \mscr{M}_{\mrm{B}} \li (
\wh{\bs{p}}_s, \udl{\bs{p}}_s \ri ) \leq \ln (\ze \de) \ri \} = 1,
\]
\[
1 \geq \Pr \li \{  S_{\mrm{B}} (0, K_s, n_s, \ovl{\bs{p}}_s ) \leq
\ze \de \ri \} \geq \Pr \li \{  n_s \mscr{M}_{\mrm{B}} \li (
\wh{\bs{p}}_s, \ovl{\bs{p}}_s \ri ) \leq \ln (\ze \de) \ri \} = 1.
\]
The lemma immediately follows.  \epf

By a similar argument as that of Lemma \ref{good1}, we have Lemma
\ref{goodmix1} as follows.

\beL \la{goodmix1} $\Pr \{ p \leq \udl{\bs{p}}_\ell, \;  S_{\mrm{B}}
(K_\ell, n_\ell, n_\ell, \udl{\bs{p}}_\ell  ) \leq \ze \de \} \leq
\ze \de$ for $\ell = 1, \cd, s$.  \eeL

By a similar argument as that of Lemma \ref{good2}, we have Lemma
\ref{goodmix2} as follows.

 \beL
 \la{goodmix2}
$\Pr \{ p \geq \ovl{\bs{p}}_\ell, \;  S_{\mrm{B}} \li (0, K_\ell,
n_\ell, \ovl{\bs{p}}_\ell \ri ) \leq \ze \de \} \leq \ze \de$ for
$\ell = 1, \cd, s$.  \eeL

\bsk

Now we are in a position to prove Theorem 4.  By the assumption that
$0 < \vep_a < \f{35}{94}$ and {\small $\f{70 \vep_a}{35 - 24 \vep_a
} < \vep_r < 1$}, we have $\f{\vep_a}{\vep_r} + \f{12}{35} \vep_a <
\f{1}{2}$. Hence, $\f{\vep_a}{\vep_r} + \vep_a < \f{1}{2} +
\f{23}{35} \vep_a < \f{1}{2} + \f{23}{35} \times \f{35}{94} < 1$. As
a result, $\vep_a + \vep_r \vep_a - \vep_r < 0$, leading to $\nu <
0$.  It follows that $\tau \leq - 1$ and thus the sample sizes $n_1,
\cd, n_s$ are well-defined.   By Lemma \ref{mixDS1}, the sampling
must stop at some stage with index $\ell \in \{1, \cd, s\}$.
Therefore, the sampling scheme is well-defined.  By Lemmas
\ref{goodmix1}, \ref{goodmix2} and the definition of the stopping
rule, we have \bee \Pr \{  | \wh{\bs{p}} - p | \geq \vep_a, \; |
\wh{\bs{p}} - p | \geq \vep_r p \} & = & \Pr \{ p \leq \udl{\bs{p}}
\} +  \Pr \{ p \geq
\ovl{\bs{p}} \}\\
& \leq & \sum_{\ell = 1}^s \Pr \{ p \leq \udl{\bs{p}}_\ell, \;
S_{\mrm{B}} (K_\ell, n_\ell, n_\ell,
\udl{\bs{p}}_\ell  ) \leq \ze \de  \} \\
&  & +  \sum_{\ell = 1}^s \Pr \li \{ p \geq \ovl{\bs{p}}_\ell, \;
S_{\mrm{B}} \li (0, K_\ell, n_\ell, \ovl{\bs{p}}_\ell \ri )
\leq \ze \de \ri \}\\
& \leq &  s \ze \de + s \ze \de = 2 s \ze \de \leq 2 (1 - \tau )
\de, \eee from which it can be seen that $\Pr \{ | \wh{\bs{p}} - p |
< \vep_a \; \tx{or} \; | \wh{\bs{p}} - p | < \vep_r p \} > 1 - \de$
if $0 < \ze < \f{1}{2 (1 - \tau)}$.  This concludes the proof of
Theorem 4.

\sect{Proof of Theorem 5}

To prove the theorem, we need to introduce function
\[ \mscr{M}_{\mrm{P}} (z, \lm) = \bec z - \lm + z \ln \f{\lm}{z}  & \tx{for} \; z > 0 \; \tx{and} \; \lm > 0,\\
- \lm & \tx{for} \; z = 0 \; \tx{and} \; \lm > 0,\\
- \iy & \tx{for} \; z \geq 0 \; \tx{and} \; \lm \leq 0.  \eec
\]

We need some preliminary results as follows.  The following results,
stated as Lemma \ref{lem99}, has been established by Chen in
\cite{Chen_EST}.

\beL \la{lem99} {\small $S_{\mrm{P}} (0, k, n, \lm) \leq \exp ( n
\mscr{M}_{\mrm{P}} ( \f{k}{n}, \lm ) )$} for $ 0 \leq k \leq n \lm$.
Similarly, {\small $S_{\mrm{P}} (k, \iy, n, \lm) \leq \exp  ( n
\mscr{M}_{\mrm{P}} ( \f{k}{n}, \lm ) )$} for $k \geq n \lm$. \eeL

\beL Let $K = \sum_{i=1}^n X_i$ where $X_1, \cd, X_n$ are i.i.d.
Poisson random variables with mean $\lm$.  Then, $\Pr \li \{
S_{\mrm{P}} \li (0, K, n, \lm \ri ) \leq \al \ri \} \leq \al$ for
any $\al > 0$. \eeL

\bpf If $\{S_{\mrm{P}}  (0, K, n, \lm ) \leq \al \}$ is an
impossible event, then $\Pr  \{ S_{\mrm{P}} (0, K, n, \lm ) \leq \al
\} = 0 < \al$.  Otherwise, if $ \{S_{\mrm{P}} \li (0, K, n, \lm \ri
) \leq \al \}$ is a possible event, then there exists an integer
$k^* = \max \{k : k  \geq 0, \; S_{\mrm{P}}  (0, k, n, \lm ) \leq
\al \}$ and it follows that $ \Pr \{ S_{\mrm{P}}  ( 0, K, n, \lm )
\leq \al \} = S_{\mrm{P}} ( 0, k^*, n, \lm  ) \leq \al$. The proof
is thus completed. \epf

\beL Let $K = \sum_{i=1}^n X_i$ where $X_1, \cd, X_n$ are i.i.d.
Poisson random variables with mean $\lm$.  Then, $\Pr \li \{
S_{\mrm{P}} \li ( K, \iy, n, \lm \ri ) \leq \al \ri \} \leq \al$ for
any $\al > 0$.

\eeL

\bpf Since {\small $\{S_{\mrm{P}} ( K, \iy, n, \lm ) \leq \al \}$}
is a possible event for any $\al > 0$, there exists an integer
{\small $k_\star = \min \{k : k \geq 0, \; S_{\mrm{P}} (k, \iy, n,
\lm  ) \leq \al \}$} and it follows that {\small $\Pr \{ S_{\mrm{P}}
( K, \iy, n, \lm ) \leq \al \} = S_{\mrm{P}} ( k_\star, \iy, n, \lm
) \leq \al$}. The proof is thus completed. \epf

The following result, stated as Lemma \ref{lem88p}, has been
established by Chen in \cite{Chen_EST}.

\beL \la{lem88p} $\Pr \li \{ \mscr{M}_{\mrm{P}} \li (
\wh{\bs{\lm}}_s, \udl{\bs{\lm}}_s \ri ) \leq \f{\ln (\ze \de) } {
n_s }, \; \mscr{M}_{\mrm{P}} \li ( \wh{\bs{\lm}}_s, \ovl{\bs{\lm}}_s
\ri ) \leq \f{\ln (\ze \de) } { n_s }  \ri \} = 1$.  \eeL

\beL \la{PoDS1} $\Pr \li \{  S_{\mrm{P}} (0, K_s - 1, n_s,
\udl{\bs{\lm}}_s ) \geq 1 - \ze \de \ri \} = \Pr \li \{ S_{\mrm{P}}
(0, K_s, n_s, \ovl{\bs{\lm}}_s ) \leq \ze \de \ri \} = 1$. \eeL

\bpf

By Lemmas \ref{lem99} and \ref{lem88p}, {\small \[ 1 \geq \Pr \li \{
S_{\mrm{P}} (0, K_s - 1, n_s, \udl{\bs{\lm}}_s ) \geq 1 - \ze \de
\ri \} = \Pr \li \{  S_{\mrm{P}} (K_s, \iy, n_s, \udl{\bs{\lm}}_s )
\leq \ze \de \ri \} \geq \Pr \li \{ \mscr{M}_{\mrm{P}} \li (
\wh{\bs{\lm}}_s, \udl{\bs{\lm}}_s \ri ) \leq \f{\ln (\ze \de)}{n_s}
\ri \} = 1,
\]}
\[
1 \geq  \Pr \li \{  S_{\mrm{P}} (0, K_s, n_s, \ovl{\bs{\lm}}_s )
\leq \ze \de \ri \} \geq \Pr \li \{ \mscr{M}_{\mrm{P}} \li (
\wh{\bs{\lm}}_s, \ovl{\bs{\lm}}_s \ri ) \leq \f{\ln (\ze \de)}{n_s}
\ri \} = 1.
\]
The lemma immediately follows.  \epf

By a similar argument as that of Lemma \ref{good1}, we have Lemma
\ref{goodpo1} as follows.

\beL \la{goodpo1} $\Pr \{ \lm \leq \udl{\bs{\lm}}_\ell, \;
S_{\mrm{P}} \li (K_\ell, \iy, n_\ell, \udl{\bs{\lm}}_\ell \ri ) \leq
\ze \de \} \leq \ze \de$ for $\ell = 1, \cd, s$.  \eeL

By a similar argument as that of Lemma \ref{good2}, we have Lemma
\ref{goodpo2} as follows.

 \beL
 \la{goodpo2}
$\Pr \{ \lm \geq \ovl{\bs{\lm}}_\ell, \;  S_{\mrm{P}} \li (0,
K_\ell, n_\ell, \ovl{\bs{\lm}}_\ell \ri ) \leq \ze \de \} \leq \ze
\de$ for $\ell = 1, \cd, s$.  \eeL

\bsk

Now we are in a position to prove Theorem 5.  Since $\f{\vep_r}{1 +
\vep_r} < \ln (1 + \vep_r) < \vep_r$ for $\vep_r \in (0, 1)$, we
have $0 < (1 + \vep_r) \ln (1 + \vep_r) - \vep_r < \vep_r^2$ and
thus $\nu > \f{1}{\vep_a \vep_r} > 1$ for $\vep_a, \; \vep_r \in (0,
1)$.  It follows that $\tau \geq  1$ and thus the sample sizes $n_1,
\cd, n_s$ are well-defined.   By Lemma \ref{PoDS1}, the sampling
must stop at some stage with index $\ell \in \{1, \cd, s\}$.
Therefore, the sampling scheme is well-defined. By Lemmas
\ref{goodpo1}, \ref{goodpo2} and the definition of the stopping
rule, we have \bee \Pr \{  | \wh{\bs{\lm}} - \lm | \geq \vep_a, \; |
\wh{\bs{\lm}} - \lm | \geq \vep_r \lm \} & = & \Pr \{ \lm \leq
\udl{\bs{\lm}} \} +  \Pr \{ \lm \geq
\ovl{\bs{\lm}} \}\\
& \leq & \sum_{\ell = 1}^s \Pr \li \{ \lm \leq \udl{\bs{\lm}}_\ell,
\; S_{\mrm{B}} \li (0, K_\ell - 1, n_\ell,
\udl{\bs{\lm}}_\ell \ri ) \geq 1 - \ze \de \ri \} \\
&  & +  \sum_{\ell = 1}^s \Pr \li \{ \lm \geq \ovl{\bs{\lm}}_\ell,
\; S_{\mrm{B}} \li (0, K_\ell, n_\ell, \ovl{\bs{\lm}}_\ell \ri )
\leq \ze \de \ri \}\\
& \leq &  s \ze \de + s \ze \de = 2 s \ze \de \leq 2 (\tau + 1 )
\de, \eee from which it can be seen that $\Pr \{ | \wh{\bs{\lm}} -
\lm | < \vep_a \; \tx{or} \; | \wh{\bs{\lm}} - \lm | < \vep_r \lm \}
> 1 - \de$ if $0 < \ze < \f{1}{2 (\tau + 1) }$.  This concludes the proof
of Theorem 5.

\end{document}